\documentclass[preprint,review,10pt]{elsarticle}
\usepackage{amsfonts}
\usepackage{amsmath,amssymb}
\usepackage{graphicx}
\usepackage{setspace}
\usepackage[left=2cm,top=2cm,right=2cm]{geometry}
\usepackage{subfig}
\newtheorem{corollary}{Corollary}

\newtheorem{lemma}{Lemma}

\newtheorem{remark}{Remark}

\newtheorem{theorem}{Theorem}
\numberwithin{equation}{section}
\journal{XYZ}
\begin{document}
\title{Approximation properties of modified Sz\'{a}sz-Mirakyan operators in polynomial weighted  space}
\author[label1,label2,slabel*]{Prashantkumar Patel}
\ead{prashant225@gmail.com}
\author[label1,label3]{Vishnu Narayan Mishra}
\ead{vishnu\_narayanmishra@yahoo.co.in;
vishnunarayanmishra@gmail.com}
\author[label4]{Mediha \"{O}rkc\"{u}}
\ead{medihaakcay@gazi.edu.tr}
\address[label1]{Department of Applied Mathematics \& Humanities,
S. V. National Institute of Technology, Surat-395 007 (Gujarat),
India}
\address[label2]{Department of Mathematics, St. Xavier's College(Autonomous), Ahmedabad-380 009 (Gujarat), India}
\address[label4]{Department of Mathematics, Faculty of Sciences, Gazi University, Teknikokullar, 06500 Ankara, Turkey}
\address[label3]{L. 1627 Awadh Puri Colony Beniganj, Phase-III, Opposite - Industrial Training Institute (ITI), Ayodhya Main Road, Faizabad, Uttar Pradesh 224 001, India}
\fntext[label*]{Corresponding authors}
\begin{abstract}
We introduce certain modified Sz\'{a}sz-Mirakyan operators in
polynomial weighted spaces of functions of one variable. We
studied approximation properties of these operators.
\end{abstract}
\begin{keyword}Sz\'{a}sz-Mirakyan operators; Rate of convergence; Weighted approximation; Polynomial Weight  \\
\textit{2000 Mathematics Subject Classification: } primary 41A25,
41A30, 41A36. \end{keyword} \maketitle
\section{Introduction}
\label{intro} Let $J_n^{[\beta]}$ be the Jain operators
\begin{equation}\label{16.1.eq1}
J_n^{[\beta]}(f,x) =  \sum_{k=0}^{\infty} \omega_{\beta}(k,n x)
f\left(\frac{k}{n}\right),
\end{equation}
where $x\in \mathbf{R}_0:=[0,\infty)$, $n\in \mathbf{N}$, $0 \leq
\beta<1$ and
\begin{equation}\label{16.1.eq3}
\omega_{\beta}(k,\alpha) = \frac{\alpha}{k!}(\alpha+
k\beta)^{k-1}e^{-(\alpha+k\beta)}  \textrm{ for } \alpha \in
\mathbf{R}_0, k\in \mathbf{N}_0=N\cup \{0\}.
\end{equation}
Approximation properties of $J_n^{[\beta]}$ were examined by Jain
\cite{jain1972approximation} for $f\in
C\left(\mathbf{R}_0\right)$. In the particular case $\beta=0$,
$J_n^{[\beta]}$ turn out to well known the Sz\'{a}sz-Mirakyan
operators \cite{szasz1950generalization}. Kantorovich type
extension of the operators (\ref{16.1.eq1}) was discussed in
\cite{umar1985approximation}. Various other generalization and its
approximation properties of similar type of operators are studied
in
\cite{agratini2013approximation,agratini2014approximation,mishrasome2013,orkcu2013szasz,Patelmishra20131,patel2015On,rempulska2009approximation,tarabie2012jain,bardaro2006approximation,bardaro2009voronovskaya}.
In this paper, we modify operators $J_n^{[\beta]}$ given by
(\ref{16.1.eq1}), i.e. we consider operators
\begin{equation}\label{16.1.eq4}
J_n^{[\beta]}(f; a_n,b_n ; x) =  \sum_{k=0}^{\infty}
\omega_{\beta}(k,a_n x) f\left(\frac{k}{b_n}\right), x\in
\mathbf{R}_0, n\in \mathbf{N}
\end{equation}
for $f\in C\left([0,\infty)\right)$, where $(a_n)_{n=1}^{\infty}$
and $\left(b_n\right)_{n=1}^{\infty}$ are given increasing and
unbounded numerical sequence such that $a_n\geq 1$, $b_n\geq 1$
and $\displaystyle \left(\frac{a_n}{b_n}\right)_1^{\infty}$ is non
decreasing and
\begin{equation}\label{16.1.eq5.1}
\frac{a_n}{b_n} = 1+ o\left(\frac{1}{b_n}\right).
\end{equation}
If $a_n = b_n = n$ for all $n \in \mathbf{N}$, then the operators (\ref{16.1.eq4}) reduce to the operators (\ref{16.1.eq1}).\\
\indent The paper is organized as follows. In our manuscript, we
shall study approximation properties of operators
(\ref{16.1.eq4}). In section 2, we shall examine moments of the
operators $J_n^{[\beta]}(f; a_n,b_n ; x)$. We discuss
approximation properties of the operators (\ref{16.1.eq4}) in
section 3. We mention Kantorovich type extension of the operators
$J_n^{[\beta]}(f; a_n,b_n ; x)$ for further research.
\section{Moments of $J_n^{[\beta]}(f;a_n,b_n;x)$}
In order to obtain moments of $J_n^{[\beta]}(f;a_n,b_n;x)$, we
need some background results, which are as follows:
\begin{lemma}[\cite{jain1972approximation}]
Let $0< \alpha <\infty$, $0\leq \beta<1$ and let the generalized
poisson distribution given by (\ref{16.1.eq3}). Then
\begin{equation}\label{16.1.eq3.1}
\sum_{k=0}^{\infty} \omega_{\beta}(\alpha,k) = 1.
\end{equation}
\end{lemma}
\begin{lemma}[\cite{jain1972approximation}]
Let $0< \alpha <\infty$, $0\leq \beta<1$. Suppose that
$$S(r,\alpha,\beta):=\sum_{k=0}^{\infty} (\alpha+\beta k)^{k+r-1} \frac{e^{-(\alpha+\beta k)}}{k!}, r=0,1,2,\ldots$$
and
$$\alpha S(0,\alpha,\beta):=1.$$
Then
\begin{equation}\label{16.1.eq3.2}S(r,\alpha,\beta) = \alpha S(r-1,\alpha,\beta) + \beta S(r,\alpha+\beta,\beta).\end{equation}
Also,
\begin{equation}\label{16.1.eq3.3}S(r,\alpha,\beta) = \sum_{k=0}^{\infty}\beta^k (\alpha+k\beta) S(r-1,\alpha+k\beta,\beta).\end{equation}
\end{lemma}
From (\ref{16.1.eq3.2}) and (\ref{16.1.eq3.3}), when $0\leq
\beta<1$, we get
\begin{eqnarray} S(1,\alpha,\beta) &=&\frac{1}{1-\beta};\label{16.1.eq3.4}\\
S(2,\alpha,\beta) &=&\frac{\alpha}{(1-\beta)^2}+\frac{\beta^2}{(1-\beta)^3};\label{16.1.eq3.5}\\
S(3,\alpha,\beta) &=& \frac{\alpha^2}{(1-\beta)^3}+\frac{\alpha\beta^2}{(1-\beta)^4}+\frac{\beta^3+2\beta^4}{(1-\beta)^5};\label{16.1.eq3.6}\\
S(4,\alpha,\beta) &=&
\frac{\alpha^3}{(1-\beta)^4}+\frac{6\alpha^2\beta^2}{(1-\beta)^5}+\frac{4\alpha\beta^3+11\alpha\beta^4}{(1-\beta)^6}+
\frac{\beta^4+ 8\beta^5+ 6\beta^6
}{(1-\beta)^7}.\label{16.1.eq3.7}
\end{eqnarray}
In the following lemma, we have computed moments up to
$4^{\textrm{th}}$ order.
\begin{lemma}\label{16.1.lemma3}
Let  $0 \leq \beta<1$, then the following equalities hold:
\begin{enumerate}
\item $\displaystyle J_n^{[\beta]}(1; a_n,b_n ; x)=1;$
\item $\displaystyle J_n^{[\beta]}(t; a_n,b_n ; x)= \frac{a_n x}{b_n(1-\beta)};$
\item $\displaystyle  J_n^{[\beta]}(t^2; a_n,b_n ; x) = \frac{x^2 a_n^2}{(1-\beta )^2 b_n^2} +\frac{x a_n}{(1-\beta )^3 b_n^2};$
\item $\displaystyle J_n^{[\beta]}(t^3; a_n,b_n ; x) = \frac{x^3 a_n^3}{(1-\beta )^3 b_n^3}+\frac{3 x^2 a_n^2}{(1-\beta )^4 b_n^3} +\frac{x (1+2 \beta ) a_n}{(1-\beta )^5 b_n^3};$
\item $\displaystyle  J_n^{[\beta]}(t^4; a_n,b_n ; x) = \frac{x^4 a_n^4}{(1-\beta )^4 b_n^4}+ \frac{6 x^3 a_n^3}{(1-\beta )^5 b_n^4}+\frac{x^2 (7+8 \beta ) a_n^2}{(1-\beta )^6 b_n^4} +\frac{x \left(1+8 \beta +6 \beta ^2\right) a_n}{(1-\beta )^7 b_n^4}.$
\end{enumerate}
\end{lemma}
\textbf{Proof: } Using equalities (\ref{16.1.eq3.1}),
(\ref{16.1.eq3.4})  to (\ref{16.1.eq3.7}) and by simple
commutation, we obtain
\begin{eqnarray*}
J_n^{[\beta]}(1; a_n,b_n ; x)& = & \sum_{k=0}^{\infty} \omega_{\beta}(k,a_n x) = 1;\\
J_n^{[\beta]}(t; a_n,b_n ; x)   &=&\frac{a_n x}{b_n}\sum_{k=0}^{\infty} \frac{1}{k!}(a_n x+ k\beta+\beta)^{k}e^{-(a_n x+k\beta+\beta)}\\
&=&\frac{a_n x}{b_n} S(1,a_n x +\beta, \beta)\\
&=& \frac{a_n x}{b_n(1-\beta)};\\
J_n^{[\beta]}(t^2; a_n,b_n ; x)
&=&\sum_{k=0}^{\infty} \frac{a_n x}{k!}(a_n x+ k\beta)^{k-1}e^{-(a_n x+k\beta)}   \frac{k^2}{b_n^2}\\
&=&\frac{a_n x}{b_n^2}\left[ S(1,a_n x +\beta, \beta) + S(2, a_n x + 2\beta,\beta)\right]\\
&=& \frac{a_n x}{b_n^2}\left[ \frac{1}{1-\beta} + \frac{a_n x + 2\beta}{(1-\beta)^2}+ \frac{\beta^2}{(1-\beta)^3}\right]\\
&=& \frac{x^2 a_n^2}{(1-\beta )^2 b_n^2} +\frac{x a_n}{(1-\beta )^3 b_n^2};\\
J_n^{[\beta]}(t^3; a_n,b_n ; x) &=&\sum_{k=0}^{\infty} \frac{a_n x}{k!}(a_n x+ k\beta)^{k-1}e^{-(a_n x+k\beta)}   \frac{k^3}{b_n^3}\\
&=& \frac{a_n x}{b_n^3}\left[ S(1,a_n x +\beta, \beta) + 3 S(2,a_n x + 2 \beta, \beta) + S(3,a_n x +3 \beta, \beta)\right]\\
&=& \frac{x^3 a_n^3}{(1-\beta )^3 b_n^3}+\frac{3 x^2 a_n^2}{(1-\beta )^4 b_n^3} +\frac{x (1+2 \beta ) a_n}{(1-\beta )^5 b_n^3};\\
J_n^{[\beta]}(t^4; a_n,b_n ; x) &=&\sum_{k=0}^{\infty} \frac{a_n x}{k!}(a_n x+ k\beta)^{k-1}e^{-(a_n x+k\beta)}   \frac{k^4}{b_n^4}\\
&=& \frac{a_n x}{b_n^4}\left[S(1,a_n x +\beta, \beta) + 7 S(2,a_n x +2\beta, \beta)\right.\\
 &&\left.+ 6S(3,a_n x +3 \beta, \beta) + S(4,a_n x +4\beta, \beta) \right]\\
&=& \frac{x^4 a_n^4}{(1-\beta )^4 b_n^4}+ \frac{6 x^3
a_n^3}{(1-\beta )^5 b_n^4}+\frac{x^2 (7+8 \beta ) a_n^2}{(1-\beta
)^6 b_n^4} +\frac{x \left(1+8 \beta +6 \beta ^2\right)
a_n}{(1-\beta )^7 b_n^4}.
\end{eqnarray*}
\begin{lemma}
Let  $0 \leq \beta<1$, then the following equalities hold:
\begin{enumerate}
\item $\displaystyle J_n^{[\beta]}(t-x; a_n,b_n ; x)= \left(\frac{a_n }{b_n(1-\beta)}-1\right)x;$
\item $\displaystyle  J_n^{[\beta]}((t-x)^2; a_n,b_n ; x)=x^2\left(\frac{a_n}{(1-\beta ) b_n}-1\right)^2+\frac{x a_n}{(1-\beta )^3 b_n^2};$
\item $\displaystyle J_n^{[\beta]}((t-x)^3; a_n,b_n ; x) = x^3 \left(\frac{a_n}{(1-\beta ) b_n}-1\right)^3+\frac{3x^2a_n}{b_n^2(1-\beta)^3} \left(\frac{a_n}{(1-\beta ) b_n}-1\right)\displaystyle + \frac{xa_n(1+2\beta)}{(1-\beta )^5 b_n^3};$
\item $\displaystyle  J_n^{[\beta]}((t-x)^4; a_n,b_n ; x) =   x^4\left(\frac{a_n}{(1-\beta )b_n}-1\right)^4 + \frac{6a_nx^3}{(1-\beta)^3b_n^2}\left(\frac{a_n}{(1-\beta )b_n}-1\right)^2\\
  ~~~~~~~~~~~~~~~~~~~~~~~~~~~~~~~~+ \frac{a_nx^2}{(1-\beta)^5b_n^3}\left(\frac{a_n(7+8\beta)}{(1-\beta)b_n}-4-8\beta\right)+x \left(\frac{a_n(1+ 8 \beta  + 6 \beta ^2)}{(1-\beta )^7 b_n^4}\right).$
\end{enumerate}
\end{lemma}
Proof of the above lemma, follows from the linearity of the
operators
$J_n^{[\beta]}(f; a_n,b_n ; x)$.\\
   By equality (\ref{16.1.eq5.1}) and $\displaystyle \lim_{n\to \infty}\beta_n =0$, we obtain
    \begin{eqnarray*}
\lim_{n\to\infty} b_n J_n^{[\beta_n]}(t-x; a_n,b_n ; x) &=& 0;\\
\lim_{n\to\infty} b_n  J_n^{[\beta_n]}((t-x)^2; a_n,b_n ; x) &=&  x;\\
\lim_{n\to\infty} b_n  J_n^{[\beta_n]}((t-x)^3; a_n,b_n ; x) &=&  0;\\
\lim_{n\to\infty} b_n^2  J_n^{[\beta_n]}((t-x)^4; a_n,b_n ; x) &=&  3x^2,\\
    \end{eqnarray*}
    for every $x\in \mathbf{R}_0$.
\section{Approximation properties}
Let $p\in \mathbf{N}_{0},$%
$$
\omega _{0}\left( x\right) =1,~~~~~~\omega _{p}\left( x\right)
=\left(1+x^{p}\right) ^{-1}~~ if~~p\geq 1,
$$
for $x\in \mathbf{R}_{0}$, and $B_{p}$ be the set of all functions $f:%
\mathbf{R}_{0}\rightarrow \mathbf{R}$ for which $f\omega _{p}$ is
bounded on $\mathbf{R}_{0}$ and the norm is given by the following
formula:
$$
\left\Vert f\right\Vert _{p}=\sup_{x\in \mathbf{R}_{0}}\omega
_{p}\left( x\right) \left\vert f\left( x\right) \right\vert .
$$
Moreover, $C_{p}$ be the set of all $f\in B_{p}$ for which
$f\omega _{p}$ is a uniformly continuous function on
$\mathbf{R}_{0}.$ The spaces $B_{p}$ and $C_{p}$ are called
polynomial weighted spaces.
\begin{lemma}\label{16.1.lemma1}
Let $r\in \mathbf{N}$ be fixed number. Then there exists positive
numerical coefficients  $\lambda_{r,j,\beta}$, $1\leq j\leq r,$
depending only on $r$ and $j$ such that
$$ J_n^{[\beta]}(t^r; a_n,b_n ; x)= \frac{1}{b_n^r(1-\beta)^r}\sum_{j=1}^{r} \frac{\lambda_{r,j,\beta}}{(1-\beta)^{j-1}} (a_n x)^j,$$
for all $x\in \mathbf{R}_0$ and $n\in \mathbf{N}$. Moreover, we
have $\lambda_{r,1,\beta} =1 = \lambda_{r,r,\beta}$.
\end{lemma}
The proof follows by a mathematical induction argument.
\begin{lemma}
For given $p\in \mathbf{N}_0$ and $(a_n)_{n=1}^{\infty}$ and
$(b_n)_{n=1}^{\infty}$ there exists a positive constant
$M_1(b_1,p,\beta)$ such that
\begin{equation}\label{16.1.eq6}\bigg\| J_n^{[\beta]}\left(\frac{1}{\omega_p(t)}; a_n,b_n ; \cdot\right)\bigg\|_p\leq M_1(b_1,p,\beta),~~~~~n\in \mathbf{N}.\end{equation}
Moreover, for every $f\in C_p$, we have
\begin{equation}\label{16.1.eq7}\| J_n^{[\beta]}(f; a_n,b_n ; \cdot)\|_p\leq M_1(b_1,p,\beta)\|f\|_p,~~~~~n\in \mathbf{N}.\end{equation}
The formula (\ref{16.1.eq3}), (\ref{16.1.eq4}) and the inequality
(\ref{16.1.eq7}), show that $J_n^{[\beta]}$, $n\in \mathbf{N}$ is
a positive linear operator from the space $C_p$  into $C_p$, $p\in
\mathbf{N}_0$.
\end{lemma}
\textbf{Proof: }
If $p=0$, then $\displaystyle \bigg\| J_n^{[\beta]}\left(\frac{1}{\omega_0(t)}; a_n,b_n ; \cdot\right)\bigg\|_0 =\displaystyle \sup_{x\in \mathbf{R}_0} | J_n^{[\beta]}(1; a_n,b_n ; x)| =1$.\\
If $p\geq 1$, then by (\ref{16.1.eq4}), (\ref{16.1.eq5.1}), Lemma
\ref{16.1.lemma3} and Lemma \ref{16.1.lemma1}, we get
\begin{eqnarray*}
\omega_p(x) J_n^{[\beta]}\left(\frac{1}{\omega_p(t)}; a_n,b_n ; x\right) &=& \omega_p(x) \left\{1+J_n^{[\beta]}(t^p; a_n,b_n ; x)\right\}\\
 &=& \frac{1}{1+x^p}\left\{1+ \frac{1}{b_n^p(1-\beta)^p}\sum_{j=1}^{p} \frac{\lambda_{r,j,\beta}}{(1-\beta)^{j-1}} (a_n x)^j \right\}\\
 &=&  \frac{1}{1+x^p} +   \frac{1}{(1-\beta)^p}\sum_{j=1}^{p} \frac{\lambda_{r,j,\beta}}{(1-\beta)^{j-1}}\frac{1}{b_n^{p-j}} \left(\frac{a_n}{b_n}\right)^j \frac{x^j}{1+x^p}\\
 &\leq & 1 +   \frac{1}{(1-\beta)^p}\sum_{j=1}^{p} \frac{\lambda_{r,j,\beta}}{(1-\beta)^{j-1}}\frac{1}{b_1^{p-j}}= M_1(b_1,p,\beta),
\end{eqnarray*}
for all $x\in \mathbf{R}_0$ and $n\in \mathbf{N}$. From this, (\ref{16.1.eq6}) follows.\\
By (\ref{16.1.eq4}) and definition of norm, we have
$$\| J_n^{[\beta]}(f; a_n,b_n ; \cdot)\|_p\leq \| J_n^{[\beta]}(\frac{1}{\omega_p(t)}; a_n,b_n ; \cdot)\|_p\|f\|_p,$$
for every $f\in C_p$, $p\in \mathbf{N}$ and $n\in \mathbf{N}$.
From (\ref{16.1.eq6}), the inequalities (\ref{16.1.eq7}) is
achieved.
\begin{theorem}
For every $p\in \mathbf{N}_0$ there exists a positive constant
$M_2(b_1,p,\beta)$ such that
\begin{eqnarray}\label{16.1.eq8}
\omega_p(x)J_n^{[\beta]}
\left(\frac{(t-x)^2}{\omega_p(t)};a_n,b_n;x\right)&\leq&
M_2(b_1,p,\beta) \left[ x^2\left(\frac{a_n}{(1-\beta )
b_n}-1\right)^2+\frac{x}{(1-\beta )^3 b_n} \right],
\end{eqnarray}
for all $x\in \mathbf{R}_0$ and $n\in \mathbf{N}.$
\end{theorem}
\textbf{Proof: }
If $p=0$, then (\ref{16.1.eq8}) follows from values of $J_n^{[\beta]} \left((t-x)^2;a_n,b_n;x\right)$.\\
Let $J_n^{[\beta]} \left(f ;x\right)= J_n^{[\beta]}
\left(f;a_n,b_n;x\right)$. Notice that
\begin{equation}\label{eq.1.eq9}
J_n^{[\beta]} \left(\frac{(t-x)^2}{\omega_p(t)};x\right)=
J_n^{[\beta]} \left((t-x)^2;x\right)+J_n^{[\beta]}
\left(t^p(t-x)^2;x\right).
\end{equation}
For $p=1$, we get
\begin{eqnarray*}
J_n^{[\beta]} \left(\frac{(t-x)^2}{\omega_1(t)};x\right)&=& J_n^{[\beta]} \left((t-x)^2;x\right)+J_n^{[\beta]} \left(t(t-x)^2;x\right)\\
&=& J_n^{[\beta]} \left((t-x)^2;x\right) +J_n^{[\beta]} \left((t-x)^3;x\right)+xJ_n^{[\beta]} \left((t-x)^2;x\right)\\
&=& (1+x)J_n^{[\beta]} \left((t-x)^2;x\right) +J_n^{[\beta]}
\left((t-x)^3;x\right).
\end{eqnarray*}
Therefore,
\begin{eqnarray*}
(1+x)J_n^{[\beta]} \left(\frac{(t-x)^2}{\omega_1(t)};x\right)&=& x^2\left(\frac{a_n}{(1-\beta ) b_n}-1\right)^2+\frac{x a_n}{(1-\beta )^3 b_n^2} + \frac{x^3}{1+x} \left(\frac{a_n}{(1-\beta ) b_n}-1\right)^3\\
&&+\frac{3x^2a_n}{(1+x)b_n^2(1-\beta)^3} \left(\frac{a_n}{(1-\beta ) b_n}-1\right)+ \frac{xa_n(1+2\beta)}{(1+x)(1-\beta )^5 b_n^3}\\
&\leq&M_2(b_1,p,\beta)\left[ x^2\left(\frac{a_n}{(1-\beta )
b_n}-1\right)^2+\frac{x}{(1-\beta )^3 b_n} \right].
\end{eqnarray*}
If $p\geq 2$, then by Lemma \ref{16.1.lemma1}, we get
\begin{eqnarray*}
\omega_p(x)J_n^{[\beta]} \left(t^p(t-x)^2;x\right)&=& \omega_p(x)\left\{J_n^{[\beta]} \left(t^{p+2};x\right)-2xJ_n^{[\beta]} \left(t^{p+1};x\right) + x^2J_n^{[\beta]} \left(t^{p};x\right)\right\}\\
&=& \frac{x}{b_n(1-\beta)}\left\{\frac{1}{b_n^{p+1}(1-\beta)^{p+1}}\sum_{j=1}^{p+1} \frac{\lambda_{p+2,j,\beta}}{(1-\beta)^{j-1}} a_n^j \frac{x^{j-1}}{1+x^p}\right.\\
&&-\frac{2}{b_n^{p}(1-\beta)^{p}}\sum_{j=1}^{p} \frac{\lambda_{p+1,j,\beta}}{(1-\beta)^{j-1}} a_n^j \frac{x^{j}}{1+x^p}\\
&& +\left. \frac{1}{b_n^{p-1}(1-\beta)^{p-1}}\sum_{j=1}^{p-1} \frac{\lambda_{p,j,\beta}}{(1-\beta)^{j-1}} a_n^j \frac{x^{j+1}}{1+x^p}\right\} + \frac{1}{(1-\beta)^{2p+3}} \left(\frac{a_n}{b_n}\right)^{p+2} \frac{x^{p+2}}{1+x^p} \\
&&- \frac{2}{(1-\beta)^{2p+1}} \left(\frac{a_n}{b_n}\right)^{p+1} \frac{x^{p+2}}{1+x^p} +  \frac{1}{(1-\beta)^{2p-1}} \left(\frac{a_n}{b_n}\right)^p \frac{x^{p+2}}{1+x^p}\\
&=& \frac{x}{b_n(1-\beta)}\left\{\frac{1}{b_n^{p+1}(1-\beta)^{p+1}}\sum_{j=1}^{p+1} \frac{\lambda_{p+2,j,\beta}}{(1-\beta)^{j-1}} a_n^j \frac{x^{j-1}}{1+x^p}\right.\\
&&-\frac{2}{b_n^{p}(1-\beta)^{p}}\sum_{j=1}^{p} \frac{\lambda_{p+1,j,\beta}}{(1-\beta)^{j-1}} a_n^j \frac{x^{j}}{1+x^p}\\
&& +\left. \frac{1}{b_n^{p-1}(1-\beta)^{p-1}}\sum_{j=1}^{p-1} \frac{\lambda_{p,j,\beta}}{(1-\beta)^{j-1}} a_n^j \frac{x^{j+1}}{1+x^p}\right\} \\
&&+\frac{x^{p+2}}{1+x^p} \left(\frac{a_n}{b_n}\right)^p
\frac{1}{(1-\beta)^{2p-1}}
\left(\frac{a_n}{b_n(1-\beta)}-1\right)^{2}.
\end{eqnarray*}
Since $\displaystyle 0\leq \frac{a_n}{b_n} \leq 1$ for $n\in
\mathbf{N}$, $\displaystyle (1-\beta)^{-1} \leq  (1-\beta)^{-3}$,
we have
\begin{eqnarray}\label{16.1.eq10}
\omega_p(x)J_n^{[\beta]} \left(t^p(t-x)^2;x\right) &\leq &
\frac{x}{b_n(1-\beta)^3}\left\{\sum_{j=1}^{p+1}
\frac{\lambda_{p+2,j,\beta}}{b_1^{p-j+1}(1-\beta)^{p+j}}  \right.
+2\sum_{j=1}^{p} \frac{\lambda_{p+1,j,\beta}}{b_1^{p-j}(1-\beta)^{p+j-1}} \nonumber\\
&& +\left. \sum_{j=1}^{p-1} \frac{\lambda_{p,j,\beta}}{b_1^{p-j-1}(1-\beta)^{p+j-2}}  \right\} + \frac{x^2 }{(1-\beta)^{2p-1}}  \left(\frac{a_n}{b_n(1-\beta)}-1\right)^{2}.\nonumber\\
&\leq& M_2(b_1,p,\beta) \left\{ x^2
\left(\frac{a_n}{b_n(1-\beta)}-1\right)^{2} +
\frac{x}{b_n(1-\beta)^3}\right\}.
\end{eqnarray}
for $x\in \mathbf{R}_0$, $n\in \mathbf{N}$. Using
(\ref{16.1.eq10}) in (\ref{eq.1.eq9}), we obtain (\ref{16.1.eq8})
for $p\geq 2$.\\ Thus, the proof is completed.\\
 Now, we approximate
$J_n^{[\beta]} \left(f;a_n,b_n;x\right)$ using the modulus of
continuity    $\omega_1(f,C_p)$ and the modulus of smoothness
$\omega_2(f,C_p)$ of function $f\in C_p$, $p\in \mathbf{N}_0$
$$\omega_1(f,C_p,t) := \sup_{0\leq h \leq t} \|\triangle_hf(\cdot)\|_p, ~~~~~~~~~~ \omega_2(f,C_p,t) := \sup_{0\leq h \leq t} \|\triangle_h^2f(\cdot)\|_p,$$
for $t\geq 0$, where
$$\triangle_hf(x) = f(x+h)-f(x), ~~~~~~\triangle_h^2 f(x) = f(x) -2f(x+h) + f(x+2h).$$
Let
\begin{equation}\label{16.1.12}
\xi_{n,\beta}(x) = x^2 \left(\frac{a_n}{b_n(1-\beta)}-1\right)^{2}
+  \frac{x}{b_n(1-\beta)^3}, ~~x\in \mathbf{R}_0, x\in \mathbf{N}.
\end{equation}
\begin{theorem}\label{16.1.thm2.5}
Suppose that $f\in C_p^2$ with a fixed $p\in \mathbf{N}_0$. Then
there exists a positive constant $M_3(b_1,p,\beta)$ such that
\begin{eqnarray}\label{16.1.eq11}
\omega_p(x)|J_n^{[\beta]} \left(f; a_n,b_n;x\right)-f(x)|\leq
\|f'\|_p\bigg|\frac{a_n}{b_n(1-\beta)}-1\bigg|~x\nonumber+
\|f''\|_p M_3(b_1,p,\beta)\xi_{n,\beta}(x),
\end{eqnarray}
for all $x\in \mathbf{R}_0$, $n\in \mathbf{N}$.
\end{theorem}
\textbf{Proof: }
Notice that $\displaystyle J_n^{[\beta]} \left(0; a_n,b_n;x\right)=f(0), n\in \mathbf{N}$, which implies (\ref{16.1.eq11}) for $x=0$.\\
Let $x>0$ and let $J_n^{[\beta]}(f;x) = J_n^{[\beta]} \left(f;
a_n,b_n;x\right)$. For $f\in C_p^2$ and $t\in \mathbf{R}_0$,
\begin{equation}
f(t) =f(x) + f'(x) (t-x) + \int_{x}^t (t-u)f''(u)du.
\end{equation}
Applying $J_n^{[\beta]}(f;x)$ on both side, we obtain
\begin{eqnarray*}
J_n^{[\beta]}(f(t);x)  =f(x) + f'(x) J_n^{[\beta]}((t-x);x) +
J_n^{[\beta]}\left(\int_{x}^t (t-u)f''(u)du;x\right).
\end{eqnarray*}
Notice that
$$\bigg| \int_{x}^t (t-u)f''(u)du\bigg| \leq \|f''\|_p \left(\frac{1}{\omega_p(t)}+\frac{1}{\omega_p(x)}\right)(t-x)^2.$$
Now, using above inequality, we have
\begin{eqnarray*}
\omega_p(x)| J_n^{[\beta]}(f(t);x)-f(x) | & \leq& \|f'\|_p J_n^{[\beta]}((t-x);x)\\
 &&+ \|f''\|_p \omega_p(x)J_n^{[\beta]}\left( \left(\frac{1}{\omega_p(t)}+\frac{1}{\omega_p(x)}\right)(t-x)^2;x\right)\\
&\leq & \|f'\|_p J_n^{[\beta]}((t-x);x)\\
 &&+ \|f''\|_p \left(\omega_p(x)J_n^{[\beta]}\left( \frac{(t-x)^2}{\omega_p(t)};x\right)+ J_n^{[\beta]}\left((t-x)^2;x\right)\right).
\end{eqnarray*}
Now, using (\ref{16.1.eq8}) and (\ref{16.1.12}), we get
\begin{eqnarray*}
\omega_p(x)| J_n^{[\beta]}(f(t);x)-f(x) | & \leq& \|f'\|_p
\bigg|\frac{a_n}{b_n(1-\beta)}-1 \bigg|~x+ \|f''\|_p
\xi_{n,\beta}(x)M_3(b_1,n,\beta).
\end{eqnarray*}
Thus, the proof is completed.
\begin{corollary}
Let $\rho(x) = (1+x^2)^{-1}$, $x\in \mathbf{R}_0$. Suppose that
$f\in C_p^2$ with a fixed $p=2$. Then there exists a positive
constant $M_4(b_1,p,\beta)$ such that
\begin{eqnarray}
\|[J_n^{[\beta]}(f;a_n,b_n;x)-f(x)]~{\rho}\|_2&\leq& \left(1-\frac{a_n}{b_n(1-\beta)}\right)\|f'\|_2\\
&&+M_4(b_1,p,\beta)\|f''\|_2b_n^{-1}(1-\beta)^{-3}, n\in
\mathbf{N}\nonumber
\end{eqnarray}
\end{corollary}
\begin{theorem}\label{16.1.thm3}
Suppose that $f\in C_p$ with a fixed $p\in \mathbf{N}_0$. Then
there exists a positive constant $M_5(b_1,p,\beta)$ such that
\begin{eqnarray*}
\omega_p| J_n^{[\beta]}(f;a_n,b_n;x) - f(x) | &\leq& \bigg|\frac{a_n}{b_n(1-\beta)}-1\bigg|x\left(\xi_{n,\beta}(x)\right)^{-1/2} \omega_1\left(f;C_p;\sqrt{\xi_{n,\beta}(x)}\right)\\
&&+
M_5(b_1,p,\beta)\omega_2\left(f;C_p;\sqrt{\xi_{n,\beta}(x)}\right),
\end{eqnarray*}
for all $x>0$ and $n\in \mathbf{N}$, where $\xi_{n,\beta}(\cdot)$
is defined in (\ref{16.1.12}). For $x=0$, it follows that
$J_n^{[\beta]}(f;a_n,b_n;0) =f(0)$.
\end{theorem}
\textbf{Proof: } We shall apply the Steklov function $f_h$ for
$f\in C_p$:
$$
f_h(x) = \frac{4}{h^2}\int_0^{h/2} \int_0^{h/2} \left[f(x+s+t) -
f(x + 2(s+t))\right] d s dt,
$$
$x\in \mathbf{R}_0$, $h>0$, for which we have
\begin{eqnarray*}
f_h'(x) &=& \frac{1}{h^2} \int_0^{h/2} \left[8\triangle_{h/2}f(x+s) - 2\triangle_hf(x + 2s)\right] d s,\\
f_h''(x) &=& \frac{1}{h^2}  \left[8\triangle_{h/2}^2f(x) -
\triangle_h^2f(x)\right].
\end{eqnarray*}
Hence, for $h>0$, we have
\begin{eqnarray}
\|f_h - f\|_p &\leq& \omega_2 (f,C_p;h), \label{16.1.eq12}\\
\|f'_h\|_p&\leq & 5h^{-1} \omega_1(f,C_p;h) \frac{\omega_p(x)}{\omega_p(x+h)},\label{16.1.eq13}\\
\|f''_h\|_p&\leq& 9h^{-2} \omega_2(f,C_p;h),\label{16.1.eq14}
\end{eqnarray}
which show that $f_h\in C_p^2$  if $f\in C_p$. By denoting
$J_n^{[\beta]}(f;a_n,b_n;x)$  by $J_n^{[\beta]}(f;x)$ we can write
\begin{eqnarray*}
\omega_p(x) | J_n^{[\beta]}(f;x)-f(x) | &\leq& \omega_p(x) \left\{|J_n^{[\beta]}(f-f_h;x)|+ |J_n^{[\beta]}(f_h;x)-f_h(x)|\right.\\
&&\left.+|f_h(x)-f(x)|\right\}:=A_1+A_2+A_3,
\end{eqnarray*}
for $x>0$, $h>0$ and $n\in \mathbf{N}$. By (\ref{16.1.eq7}) and
(\ref{16.1.eq12}), we have
\begin{eqnarray*}
A_1&\leq& M_1(b_1,p,\beta)\|f-f_h\|_p \leq M_1(b_1,p,\beta)\omega_2 (f,C_p;h),\\
A_3&\leq & \omega_2 (f,C_p;h).
\end{eqnarray*}
Applying Theorem \ref{16.1.thm2.5}, inequalities (\ref{16.1.eq13})
and (\ref{16.1.eq14}), we get
\begin{eqnarray*}
A_2&\leq&  \|f'\|_p\bigg|\frac{a_n}{b_n(1-\beta)}-1\bigg|x+ \|f''\|_p M_3(b_1,p,\beta)\xi_{n,\beta}(x)\\
&\leq & \omega_1(f,C_p;h)
\frac{\omega_p(x)}{\omega_p(x+h)}\frac{5x}{h}\bigg|\frac{a_n}{b_n(1-\beta)}-1\bigg|+
\frac{9}{h^2} \omega_2(f,C_p;h) M_3(b_1,p,\beta)\xi_{n,\beta}(x)
\end{eqnarray*}
Combining these and setting  $h=\sqrt{\xi_{n,\beta}(x)},$ for
fixed $x>0$ and $n\in \mathbf{N}$, we obtain the desired result.
\begin{theorem}\label{16.1.thm4}
Let $f\in C_p$, $p\in \mathbf{N}_0$, and let $\rho(x) =
(1+x^2)^{-1}$ for $x\in \mathbf{R}_0$. Then there exists a
positive constant $M_6(b_1,p,\beta)$ such that
\begin{eqnarray*}
\big\|\left[J_n^{[\beta]}(f;a_n,b_n;x)-f\right]   \rho\big\|_p &\leq& \left(1-\frac{a_n}{b_n(1-\beta)}\right) \sqrt{b_n} ~~ \omega_1\left(f,C_p;1/\sqrt{b_n(1-\beta)^3}\right)\\
&&+ M_6(b_1,p,\beta)
\omega_2\left(f,C_p;1/\sqrt{b_n(1-\beta)^3}\right), ~~~n\in
\mathbf{N}.
\end{eqnarray*}
\end{theorem}
From Theorems \ref{16.1.thm3} and \ref{16.1.thm4}, we derive the
following corollary:
\begin{corollary}
Let $f\in C_p$, $p\in \mathbf{N}_0$, $\beta_n\to 0$ as $n\to
\infty$. Then for $J_n^{[\beta_n]}$ defined by (\ref{16.1.eq4}),
we have
\begin{equation}\label{16.1.eq15}
\lim_{n\to \infty} J_n^{[\beta_n]} (f; a_n,b_n;x) =f(x), ~~~~x\in
\mathbf{R}_0.
\end{equation}
Furthermore, the convergence of (\ref{16.1.eq15}) is uniformly on
every interval $[x_1,x_2]$, where $ x_2>x_1\geq0$.
\end{corollary}
\begin{remark}
 The error of approximation of a function $f\in
C_{p},p\in \mathbf{N}_{0}$ by $J_{n}^{\left[ \beta \right] }\left(
f;a_{n},b_{n};.\right) $ where $a_{n}=n^{r}+\frac{1}{n}$ and
$b_{n}=n^{r},r>1 $ is smaller than by the operators
\eqref{16.1.eq1}.
\end{remark}
\section{The Operators $K_n^{[\beta]}(f;a_n,b_n)$}
In 1985, Umar and Razi \cite{umar1985approximation} introduced
Kantorovich type extension of the operators (\ref{16.1.eq1}).
Motivated by these, we introduce the further generalization of the
operators (\ref{16.1.eq4}), in the following way:
\begin{equation}\label{16.1.eq2}
K_n^{[\beta]}(f;a_n,b_n;x) =  b_n\sum_{k=0}^{\infty}
\omega_{\beta}(k,a_n x)\int_{\displaystyle
\frac{k}{b_n}}^{\displaystyle \frac{k+1}{b_n}} f\left(t\right)dt,
\end{equation}
where $x\in \mathbf{R}_0:=[0,\infty)$, $n\in \mathbf{N}$, $0 \leq \beta<1$ and $\omega_{\beta}(k,a_n x)$ as same in (\ref{16.1.eq4}).\\
 It is obvious that, the operators $K_n^{[\beta]}(f,a_n,b_n)$, $n\in \mathbf{N}$, defined in (\ref{16.1.eq2}),
 we can consider for function $f\in C_p$, $p\in \mathbf{N}_0$. For these operators and $f\in C_p$, we can prove lemma and theorems similar to the
 operators $J_n^{[\beta]}(f,a_n,b_n)$.\\
    One can study properties of $K_n^{[\beta]}(f,a_n,b_n)$ for functions $f\in L(\mathbf{R}_0)$ as same as discussed in \cite{walczak2002approximation}.


\end{document}